\documentclass[12pt,fleqn]{article}
\usepackage{latexsym}
\usepackage{amsmath}
\usepackage{amssymb}
\usepackage{epsfig}
\usepackage{amsfonts}

\begin{document}

\newtheorem{definition}{Definition}[section]
\newtheorem{theorem}[definition]{Theorem}
\newtheorem{lemma}[definition]{Lemma}
\newtheorem{remark}[definition]{Remark}
\newtheorem{remarks}[definition]{Remarks}
\newtheorem{examples}[definition]{Examples}
\newtheorem{leeg}[definition]{}
\newtheorem{definitions}[definition]{Definitions}
\newtheorem{proposition}[definition]{Proposition}
\newtheorem{example}[definition]{Example}
\newtheorem{comments}[definition]{Some comments}
\newtheorem{corollary}[definition]{Corollary}
\def\square{\Box}
\newtheorem{observation}[definition]{Observation}
\newtheorem{observations}[definition]{Observations}
\newtheorem{defsobs}[definition]{Definitions and Observations}
\newenvironment{prf}[1]{ \trivlist
\item[\hskip \labelsep{\it
#1.\hspace*{.3em}}]}{~\hspace{\fill}~$\square$\endtrivlist}
\newenvironment{proof}{\begin{prf}{Proof}}{\end{prf}}

\title{Analytic $q$-difference equations}
\author{Marius van der Put \\
\footnotesize Department of Mathematics, University of Groningen,
 P.O.Box 800,\\
\footnotesize 9700 AV Groningen, 
The Netherlands, mvdput@math.rug.nl }
\date{}
\maketitle
\noindent

\section*{Introduction} 
A complex number $q$ with $0<|q|<1$ is fixed. By an analytic $q$-difference equation we mean 
an equation which can be represented by a matrix equation $Y(z)=A(z)Y(qz)$ where $A(z)$
is an invertible $n\times n$-matrix with coefficients in the field $K=\mathbb{C}(\{z\})$ of the
convergent Laurent series and where $Y(z)$ is a vector of size $n$. The aim of this paper
is to give an overview of our present knowledge of these equations and their solutions.
Definitions and statements are presented in detail. Examples illustrate the main results. For proofs we refer to the cited literature.
For section 1 and part of the following sections, the reference is [P-S]. For the later sections the reference is [P-R]. The last section presents unpublished results. The theory of linear differential equations over
$K$ (see [P-S.2], especially Chapter 10) has many of the features presented in this survey. 

The manuscript [R-S] 
contains a concise overview of analytic $q$-difference equations and its main purpose is
to develop  a theory of $q$-summation leading to, in our terminology, a description of a
universal difference Galois group.  As the authors of [R-S] note, this is only partially achieved
and part of their work is still conjectural. However for a certain class of $q$-difference equations,
namely those having at most two slopes and such that the slopes are integral,
they have explicit results. An important part of the extensive literature on $q$-difference equations can be found in the papers cited here.

\section{Difference equations in general}
A difference field $F$ is a field provided with an automorphism $\phi$ of infinite order. 
A scalar linear difference equation is an equation of the form 
\[\phi ^n(y)+a_{n-1}\phi ^{n-1}(y)+\cdots +a_1\phi (y) +a_0y=0\,\]
with given $a_i\in F$ and, say, $a_0\neq 0$. 
As in the case of linear differential equations, one can transform this equation into a matrix difference equation, i.e., an equation of the form $Y=A\phi (Y)$, where $A$ is a given invertible $n\times n$-matrix with coefficients in $F$ and $Y$ denotes a vector of length $n$. On the vector space $M=F^n$ one considers the operator 
$\Phi : Y\mapsto A\phi (Y)$. The bijective map $\Phi :M\rightarrow M$ is additive and 
$\Phi (f\cdot m)=\phi (f)\cdot \Phi (m)$ for $m\in M$ and $f\in F$. This leads to the following definition of a {\it difference module $M=(M,\Phi )$ over $F$}:\\
$M$ is a finite dimensional vector space over $F$ and $\Phi :M\rightarrow M$ is an additive, bijective map
satisfying $\Phi (f\cdot m)=\phi (f)\cdot \Phi (m)$ for $m\in M$ and $f\in F$. The equation, corresponding
to a difference module is $Y=\Phi (Y)$. In general this equation has few solutions in $M$ itself. Like ordinary
polynomial equations over fields, one has to construct extensions of $F$ in order to have sufficiently many
solutions. 

 We will now assume that $F$ has {\it characteristic zero} and that its field of constants 
 $C: =\{f\in F|\ \phi (f)=f\}$ is {\it algebraically closed}. A {\it Picard-Vessiot ring (or extension) $R$ for a difference module $M$}, say, represented by  the matrix equation $Y=A\phi (Y)$ is defined by:\\
 (i) $R$ is an $F$-algebra (commutative and with a $1$),\\
 (ii) $R$ is provided with an automorphism $\phi$ extending $\phi$ on $F$,\\
 (iii) $R$ has only trivial $\phi$-invariant ideals,\\  
 (iv) there exists an invertible matrix $U$ (called fundamental matrix) with coefficients in $R$ such that
  $U=A\phi (U)$,\\
 (v) $R$ is generated over $F$ by the coefficients of $U$ and $\frac{1}{\det U}$.\\ 

 Property (iii) translates in terms of $M$ into $V:=\{a\in R\otimes _FM|\ \Phi (a)=a\}$ is a vector space over
 $C$ and the natural map $R\otimes _CV\rightarrow R\otimes _FM$ is an isomorphisms. 
 
 One constructs $R$ as follows. Let $X$ denote a matrix $(X_{i,j})$ of indeterminates  and put $D:=\det X$.
 The $F$-algebra $R_0:=F[X,\frac{1}{D}]$ is provided with a $\phi$-action, extending the one on $F$, by
 the formula $(\phi X_{i,j} )=A^{-1}(X_{i,j})$. Let $I\subset R_0$ denote an ideal maximal among the ideals
 invariant under $\phi$. Then $R=R_0/I$ is a Picard-Vessiot ring for the given equation.

The basic results of difference Galois theory are:\\
(1) A Picard-Vessiot ring $R$ exists and is unique up to a, non unique, isomorphism. \\  
(2) $R$ is reduced (i.e., has no nilpotent elements).\\
(3) The set constants of the ring of total fractions of $R$ is $C$.\\ 
(4) Let $G$ be the group of the $F$-linear automorphism of $R$, commuting with $\phi$.
 The natural action of $G$ on $R\otimes _FM$ induces a faithful action of $G$ on $V$, the solution space.
 The image of $G$ in ${\rm GL}(V)$ is a linear algebraic subgroup of the latter. This makes $G$ into  a  linear algebraic group over $C$.\\ 
(5) The action of $G$ on $Spec(R)$ makes the latter into an  $G$-torsor over $F$. In other words, there exists
a finite extension $F^+\supset F$ and a $G$-equivariant isomorphism 
$F^+\otimes _CC[G]\rightarrow F^+ \otimes _FR$, where $C[G]$ is the coordinate ring of $G$.\\

Most of the notions and `operations of linear algebra', such as  morphisms, kernels, cokernels, direct sums
have an obvious equivalent for difference modules. Let $(M_1,\Phi _1),\ (M_2,\Phi _2)$ denote two difference modules. The tensor product of the two modules is defined as $M_1\otimes _FM_2$ with $\Phi$ given by
$\Phi (m_1\otimes m_2)=(\Phi _1m_1)\otimes (\Phi _2m_2)$. The internal hom of the two modules is
${\rm Hom}_F(M_1,M_2)$ with $\Phi$ defined by $(\Phi (L))(m_1)=\Phi _2^{-1}(L(\Phi _1m_1))$ for
$L\in {\rm Hom}_F(M_1,M_2)$ and $m_1\in M_1$. 

This leads to another, more abstract but very useful, formulation of the above Picard-Vessiot theory, namely that of (neutral) Tannakian categories. The category of all difference modules over $F$ is a neutral Tannakian
category. We will return to this in section 6. For a specific difference field $F$ one can 
say much more than the above formalism (analogous to the case of ordinary Galois theory for specific fields).

\section{First examples of $q$-difference equations}

This exposition is concerned with the difference field $K=\mathbb{C}(\{z\})$, i.e., the field of convergent Laurent
series over $\mathbb{C}$, provided with the automorphism $\phi$ given by $\phi (z)=qz$, where
$q$ is a fixed complex number satisfying $0<|q|<1$.   In order to define $\phi$ on the algebraic closure
$K_\infty =\cup _{n\geq 1} K_n$, with $K_n:=\mathbb{C}(\{z^{1/n}\})$, of
$K$ we choose a $\tau \in \mathbb{C}$ with $\Im (\tau )>0$ such that $q=e^{2\pi i \tau}$. Define $q^\lambda$ for $\lambda \in \mathbb{Q}$ (or any $\lambda \in \mathbb{C}$) as $e^{2\pi i\lambda \tau}$. Then the action of $\phi$ on 
 $K_\infty$ is given by $\phi (z^\lambda)=q^\lambda z^\lambda $. The action of $\phi$ on 
$\widehat{K}=\mathbb{C}((z))$, i.e., the field of the formal Laurent series,   and on its algebraic closure $\widehat{K}_\infty =\cup _{n\geq 1}\widehat{K}_n$, with $\widehat{K}_n:=\mathbb{C}((z^{1/n}))$, 
is defined in a similar way.

 Some $q$-difference rings $F=(F,\phi )$ will be considered, namely $\mathbb{C}[z,z^{-1}]$ and
 $O$, the ring of the holomorphic functions on $\mathbb{C}^*$. A difference module $M=(M,\Phi )$ over a $q$-difference ring $F$ will be a {\it free} $F$-module of finite rank provided with a bijective additive map
 $\Phi :M\rightarrow M$ such that $\Phi (f\cdot m)=\phi (f)\cdot \Phi (m)$ for $f\in F,\ m\in M$.  
 
\begin{examples} {\rm $\ $\\
 (a) $Ke, \Phi (e)=ce,\ c\in \mathbb{ C}^*$. There exists a solution ($\neq 0$)
in $K$ $\Leftrightarrow \ c\in q^{\mathbb{Z}}$. Indeed, $\Phi (z^te)=cq^tz^te$. If $c\in q^{\mathbb{Z}}$, then
the Picard-Vessiot ring is $K$ itself and the difference Galois group is $\{1\}$.

There exists a solution ($\neq 0$) in $K_n$ $\Leftrightarrow$ $c\in q^{\frac{1}{n}\mathbb{Z}}$. If 
$c\in q^{\frac{1}{n}\mathbb{Z}}$ with $n$ minimal, then $K_n$ is the Picard-Vessiot ring and the difference
Galois group is $\mu _n:=\{a\in \mathbb{C}|\ a^n=1\}$.

In the remaining case, the Picard-Vessiot ring is $K[X,X^{-1}]$ with the action of $\phi$ given by
$\phi (X)=c^{-1}X$. The difference Galois group is $\mathbb{G}_m:=\mathbb{C}^*$. It consists of the automorphisms $\sigma$ given by $\sigma (X)=aX$ (for any $a\in \mathbb{C}^*$). We will use the symbol
$e(c)$ for this $X$. It can be given an interpretation as multivalued function $z^b$ for a $b$ with $q^b=c^{-1}$.

\noindent (b) The difference module  $U_n:=Ke_1+\cdots +Ke_n$ with $\Phi$ given by the matrix
\[\left[\begin{array}{cccc}1&1& &\\ 
                            &1&1 &\\
                            &&&1\\
                            &&& 1 \end{array}\right] \]

is called {\it unipotent of length $n$}. The Picard-Vessiot ring is $K[X]$,  with $\phi (X )=X +1$.
The difference Galois group is $\mathbb{G}_a:=\mathbb{C}$. The elements $\sigma$ of this group have the form $\sigma (X )=X +a$ (any $a\in \mathbb{C}$). We will use the symbol $\ell$ for this $X$. It has an interpretation as multivalued function $\frac{\log z}{2\pi i \tau}$. For $n=2$ one easily verifies that the solution space has $\mathbb{C}$-basis $\{e_1-\ell e_2,\ e_2\}$. \hfill $\square$ }\end{examples}

A difference module $M$ over $K$ is called {\it regular singular} if $M$ has a basis $e_1,\dots ,e_m$ such that
$\mathbb{C}\{z\}e_1+\cdots +\mathbb{C}\{z\}e_m$ is invariant under $\Phi$ and  $\Phi ^{-1}$.

\begin{theorem} The following are equivalent\\  
{\rm (i)} $M$ regular singular.\\
{\rm (ii)} $M=K\otimes _{\mathbb{C}}W$, $\dim _{\mathbb{C}}W < \infty$ and 
$\Phi (f\otimes w)=\phi (f)\otimes A(w)$ for some $A\in {\rm GL}(W)$. Moreover, there is a unique $A$ such that 
 every eigenvalue $c$ satisfies $|q|<|c|\leq 1$.\\
{\rm (iii)} $M$ is obtained by $\oplus ,\otimes $ from the examples in {\rm 2.1}.
\end{theorem}

The difference module $(Ke,\ \Phi e=(-z)e)$ is the basic example of an {\it irregular singular} difference module. Its Picard-Vessiot
ring is $K[X,X^{-1}]$ with $\phi (X)=(-z)^{-1}X$. The difference Galois group is $\mathbb{G}_m$. We will
use the symbol $e(-z)$ for this $X$. It has the interpretation $ \Theta (z):=\sum _{n\in \bf Z}q^{n(n-1)/2}(-z)^n$
because of the well known formula $(-z)\Theta (qz)=\Theta (z)$. 

\section{Towards a classification of modules}
We present here  the `classics' by  G.D.~Birkhof, P.E.~Guenther, C.R.~Adams et al. and `modern' work
by J.-P.~Ramis, Ch.~Zhang, J.~Sauloy, A.~Duval, M.F.~Singer, M.~van der Put, M.~Reversat et al.,
concerning $q$-difference equations over $K$.

Let $K[\Phi ,\Phi ^{-1}]$ be the skew ring of difference operators. The elements of this ring are the finite
formal sums $\sum _{n\in \mathbb{Z}}a_n\Phi ^n$ and the multiplication is given by the rule $\Phi \cdot f
=\phi (f)\Phi$. This ring is (left and right) Euclidean. Any difference module $(M,\Phi _M)$ can be seen as
a left $K[\Phi ,\Phi ^{-1}]$-module where the action of $\Phi$ on $M$ is just $\Phi _M$. This left module is 
{\it cyclic} (i.e., generated by one element) and therefore $M\cong K[\Phi ,\Phi ^{-1}]/K[\Phi ,\Phi ^{-1}]L$ for
some $L=\Phi ^m+a_{m-1}\Phi ^{m-1}+\cdots +a_1\Phi +a_0$ with $a_0\neq 0$.

The usual discrete valuation $v$ on $K$ is given by $v(0)=+\infty$ and for $a\in K^*$, $v(a)$ is the order of 
$a$. Using the values $v(a_i)$ one defines the {\it Newton polygon} of $L$, as in the case of an ordinary
polynomial in $K[T]$. This Newton polygon depends on $M$ only. The slopes of the Newton polygon are in
$\mathbb{Q}$. A difference module $M$ is called {\it pure}  if there is only one slope. 

\begin{examples} {\rm 
$M$ is regular singular if and only if $M$ is pure of slope 0.

\noindent 
Further $(Ke,\ \Phi e=c(-z)^te)$ with $c\in \mathbb{C}^*,\ t\in \mathbb{Z}$ is pure of slope $t$.  }\hfill $\square$\end{examples}
\begin{theorem}[Adams, Birkhoff, Guenther, Ramis, Sauloy] $\ $\\
$M$ has a unique tower of submodules 
$0=M_0\subset M_1\subset \cdots \subset M_r=M$ such that every $M_i/M_{i-1}$
is pure of slope $\lambda _i$ and $\lambda _1<\cdots <\lambda _r$.
\end{theorem}
This filtration is called the {\it slope filtration} of $M$ and one defines the {\it graded module} $gr(M)$ of $M$
by $gr(M)=\oplus _i M_i/M_{i-1}$.  In proving Theorem 3.2, one observes that the above operator $L$ has a unique factorization $L_1\cdot L_2\cdots L_r$ with each $L_i\in \widehat{K}[\Phi ,\Phi ^{-1}]$ monic and having only one slope $\lambda _i$. The main step is to show that these $L_i$ are actually convergent, i.e., belong
to $K[\Phi ,\Phi ^{-1}]$. This factorization of $L$ induces the slope filtration. 

We note that for obtaining convergence, it is essential that $\lambda _1<\cdots <\lambda _r$. For any order
of the slopes $\lambda _i$ there is a similar factorization of $L$ in the ring $\widehat{K}[\Phi ,\Phi ^{-1}]$.
This has as consequence that $\widehat{K}\otimes _KM$ is in fact equal to the direct sum 
$\oplus _i(\widehat{K}\otimes _K M_i/M_{i-1})=\widehat{K}\otimes _Kgr(M)$. In the next section we will study the moduli spaces that describe the modules $M$ with a fixed graded module $gr(M)$. A difference module
$M$ over $K$ is called {\it split} if it is isomorphic to its graded module $gr(M)$.  
Here we continue with the classification of pure
modules over $K$. We note that by Theorem 3.2 any irreducible module over $K$ is pure.

\begin{definition} The module $E(cz^{t/n})$.\\{\rm 
Given are  the data $t/n,\ n\geq 1,\ (t,n)=1, \ c\in \mathbb{C}^*, |q|^{1/n}<|c|\leq 1$. They define
a difference module over $K_n$ of dimension 1, namely $(K_ne, \ \Phi e=cz^{t/n}e)$. This object,
considered as difference module over $K$, has dimension $n$ over $K$ and is called
 $E(cz^{t/n})$.}\hfill $\square$\end{definition}

\begin{theorem} $E(cz^{t/n})$ is pure, irreducible and has slope $t/n$. Further,
$E(c_1z^{t/n})\cong E(c_2z^{t/n})$ if and only if $ c_1^n=c_2^n $. Moreover,
every irreducible $M$ over $K$ is isomorphic to some  $E(cz^{t/n})$.
\end{theorem}

A difference module is called {\it indecomposable} if it is not the direct sum of two proper submodules.
We note that an indecomposable module over $K$ need not be pure.

\begin{theorem} The indecomposable pure modules over $K$ are
$E(cz^{t/n})\otimes U_m$ with $|q|^{1/n}<|c|\leq 1$. Further, the triple $(t/n, c^n,m)$ is unique.
\end{theorem}

\begin{definition} Global lattices. \\{\rm
A global lattice $\Lambda$ for a $q$-difference module $M$ over $K$ is a finitely generated
$\mathbb{C}[z,z^{-1}]$-submodule of $M$ (and hence free), invariant under $\Phi$ and $ \Phi ^{-1}$, such that
the natural map $K\otimes _{\mathbb{C}[z,z^{-1}]}\Lambda \rightarrow M$ is an 
isomorphism. }\hfill $\square$\end{definition}

We will see later that {\it any} difference module $M$ has a unique global lattice. This means that the $q$-difference equation, defined locally at $z=0$, is equivalent to an equation on all of $\mathbb{P}^1$ with 
at most poles at $z=0$ and $z=\infty$. From the explicit description of the indecomposable pure modules over $K$ it is not hard to deduce the following.

\begin{corollary} Every pure indecomposable difference module $M$ over $K$ has a unique global lattice.
This lattice, {\em denoted by} $M_{global}$, is a difference module over $\mathbb{C}[z,z^{-1}]$. The same holds for {\em split} difference modules over $K$.\\
Moreover, any morphism $f:M\rightarrow N$ between split difference modules over $K$ maps 
$M_{global}$ into $N_{global}$. \end{corollary}

Before going on, we remark that the {\it classification of the difference modules over $\widehat{K}$} is remarkably  simple. With the same methods used in the proof of Theorem 3.5, one shows that
every pure indecomposable difference module over $\widehat{K}$ has the form
$\widehat{K}\otimes _K(E(cz^{t/n})\otimes U_m)$ (again with unique $(t/n,c^n,m)$). Finally, as remarked before, any difference module over $\widehat{K}$ is a direct sum of pure modules over $\widehat{K}$.\\

It is well known that the elliptic curve $E_q:={\mathbb{C}}^*/q^{\mathbb{Z}}$, which we like to call {\it the Tate curve}, plays an important role for $q$-difference equations. With the help of 3.5, 3.6 and 3.7 one can deduce the following 
rather striking result.

\begin{theorem} There is an additive, faithful functor $V$ from the category of the split difference modules over $K$ to the category of the vector bundles on $E_q$. It has the properties:\\ 
{\rm (i)} $V$ induces a bijection between the (isomorphy classes of) indecomposable modules over
$K$ and the (isomorphy classes of) indecomposable vector bundles on $E_q$.\\
{\rm (ii)} $V$ induces a bijection between (isomorphy classes of) objects. \\ 
{\rm (iii)} $V$ respects the constructions of linear algebra, i.e., tensor products, exterior powers etc.
\end{theorem}

\begin{proof} We will {\it sketch} a proof. \\
(1). We recall that $O$ denotes the algebra of the holomorphic functions on 
${\mathbb C}^*$ and that a difference module $M$ over $O$ is a left module 
over the ring $O[\Phi ,\Phi ^{-1}]$, free of some rank $m<\infty $ over 
$O$. Further $pr:{\mathbb C}^*\rightarrow E_q:={\mathbb C}^*/q^{\mathbb Z}$
denotes the canonical map. One associates 
to $M$ the vector bundle $v(M)$ of rank $m$ on $E_q$ given by
$v(M)(U)=\{f\in O(pr ^{-1}U)\otimes _OM|\ \Phi (f)=f\}$,
where, for any open $W\subset {\mathbb C}^*$, $O(W)$ is the algebra of the
holomorphic functions on $W$.

On the other hand, let a vector bundle $\mathcal M$ of rank
$m$ on $E_q$ be given. Then ${\mathcal N}:=pr^*{\mathcal M}$ is a vector
bundle on ${\mathbb C}^*$ provided with a natural isomorphism
$\sigma _q^*{\mathcal N}\rightarrow {\mathcal N}$, where $\sigma _q$ is the
map $\sigma _q(z)=qz$. One knows that ${\mathcal N}$ is in fact a free
vector bundle of rank $m$ on ${\mathbb C}^*$. Therefore, $M$, the collection
of the global sections of $\mathcal N$, is a free $O$-module
of rank $m$ provided with an invertible action $\Phi$ satisfying
$\Phi (fm)=\phi (f)\Phi (m)$ for $f\in O$ and $m\in M$.
It is easily verified that the above describes an equivalence $v$ of tensor categories.\\ 

\noindent (2). One associates to any split difference module $M$ over $K$, its global lattice $M_{global}$
and the $q$-difference module $O\otimes _{\mathbb{C}[z,z^{-1}]}M_{global}$ over $O$. The latter
induces by (1) a vector bundle on $E_q$, which we call $V(M)$. For a  morphism $f:M\rightarrow N$ between split modules one has $f:M_{global}\rightarrow N_{global}$.Therefore $f$ induces a morphism $V(f):V(M)\rightarrow V(N)$. Thus we found the additive, faithful functor $V$. Clearly $V$ respects the constructions of linear algebra.\\

\noindent (3). As (ii) is an immediate consequence of (i), we are left with proving (i). The indecomposable module $M:=E(cz^{t/n})\otimes U_m$ is producing a vector bundle $V(M)$ of rank $nm$ and degree $tm$. 
One can show that $V(M)$ is indecomposable and that $V(M)$ is irreducible if $m=1$. Further one can verify
that non isomorphic indecomposable $M, N$ produce non isomorphism vector bundles 
$V(M), V(N)$. It is somewhat more complicated to show that every indecomposable vector bundle is isomorphic to $V(M)$ for a suitable indecomposable $M$.  This last step can be avoided by an inspection of Atiyah's paper where the classification of the indecomposable vector bundles on $E_q$ is explicitly 
given.  \end{proof}

\begin{corollary} Let $B$ be a split module over $K$, then\\
{\rm (i)} ${\rm ker}(\Phi -1,O\otimes B_{global})\cong H^0(E_q,V(B))$.\\
{\rm (ii)} ${\rm coker}(\Phi -1,O\otimes B_{global})\cong H^1(E_q,V(B))$.\\
{\rm (iii)} The two canonical maps ${\rm coker}(\Phi -1,B_{global})
\rightarrow {\rm coker}(\Phi -1,B)$ and ${\rm coker}(\Phi -1,B_{global})
\rightarrow  {\rm coker}(\Phi -1,O\otimes B_{global})$ are isomorphisms.
\end{corollary}

\begin{examples}{\rm Consider the difference module
$(B=Ke,\ \Phi e=(-z)^te)$ with $t\in \mathbb{Z}$. The line bundle
$V(B)$ is equal to $O_{E_q}(t\cdot [1])$, where $1$ denotes the neutral element of $E_q$.

For $t\geq 1$, ${\rm coker}(\Phi -1,O\otimes B_{global})=0$ and
${\rm ker}(\Phi -1,O\otimes B_{global})$ is the $t$-dimensional vector space with basis
$\{ \Theta (\zeta z)^te| \zeta ^t=1\}$.

For $t=0$, ${\rm ker}(\Phi -1,O)=\mathbb{C}1$ and
${\rm coker}(\Phi -1,O)$ has dimension 1. Explicitly, $O$ is the vector space of the everywhere
convergent Laurent series and has the formula
$(\Phi -1)(\sum _{n\in \mathbb{Z}}a_nz^n)
=\sum _{n\in \mathbb{Z}}(q^n-1)a_nz^n$. }\hfill $\square$ \end{examples}

\section{Moduli spaces for difference modules}
Fix a split difference module $S=P_1\oplus \cdots \oplus P_r$ over $K$, where $P_i$ is pure
with slope $\lambda _i$ and $\lambda _1<\cdots <\lambda _r$. The problem is to classify
the difference modules $M$ over $K$ such that $gr(M)$ is isomorphic to $S$. The collection
of the isomorphy classes is a rather ugly object due to the fact that $S$ has many automorphisms. 

A better way to formulate the problem is to consider pairs $(M,f)$ consisting
of a difference module $M$ and an isomorphism $f:gr(M)\rightarrow S$. Two pairs 
$(M_1,f_1), (M_2,f_2)$ are called {\it equivalent} if there exists an isomorphism 
$g:M_1\rightarrow M_2$ such that the induced graded isomorphism $gr(g):gr(M_1)\rightarrow
gr(M_2)$ has the property $f_1=f_2\circ gr(g)$. Let $Equiv(S)$ denote the set of equivalence classes.

This formulation allows us to define a covariant functor $\mathcal F$ from the category of finitely generated $\mathbb{C}$-algebras $R$ (i.e., $R$ is commutative and has a $1$) to the category of sets. For a finitely generated $\mathbb{C}$-algebra $R$ one considers 
$K_R:=R\otimes _{\mathbb{C}}K$ and one can define the notion of difference module over
$K_R$.  The set ${\mathcal F}(R)$ is the set of equivalence classes of pairs
$(M,f)$ with $M$ a difference module over $K_R$ and $f$ is a $K_R$-isomorphism $f:gr(M)\rightarrow K_R\otimes S$. Equivalence of two pairs is defined as above.
 We note that ${\mathcal F}(\mathbb{C})$ is precisely $Equiv(S)$. 
  
One can see the functor $\mathcal F$ as a contravariant functor on the category of affine
$\mathbb{C}$-schemes (of finite type). 

\begin{theorem} The contravariant functor $\mathcal F$ on affine $\mathbb{C}$-schemes is representable. In fact, the affine space $\mathbb{A}_{\mathbb{C}}^N$ with  
$N=\sum _{i<j}(\lambda _j-\lambda _i)\dim P_i\cdot \dim P_j$, represents $\mathcal F$.
 In other words, the covariant 
functor $\mathcal F$ is represented by a certain universal $q$-difference module $M$ over  
$\mathbb{C}[X_1,\dots ,X_N]$. Further $Equiv(S)$ identifies with $\mathbb{C}^N$.
\end{theorem}
For the case of integer slopes $\lambda _i$, the above result is announced  by 
J.-P.~ Ramis and J.~Sauloy. The general case is treated in [P-R].

One can normalize the representing space $\mathbb{A}^N_{\mathbb{C}}$ by letting
$0$ correspond to the class of the pair $(S,id_S)$. The vector space structure of 
$\mathbb{A}^N_{\mathbb{C}}$ has an interpretation for $s=2$, namely as the vector space
${\rm coker} (\Phi -1,{\rm Hom}(P_2,P_1))$. For $s>2$, the functor $\mathcal F$ and its representing space 
 still have a weaker structure, namely that of an iterated torsor. 
 
We illustrate Theorem 4.1 by the following {\it basic example}:\\
$S=P_1\oplus P_2$, where $P_1=(Ke_1,\ \Phi e_1=e_1)$ and  $P_2=(Ke_2,\  \Phi e_2=(-z)^te_2)$ with  $t>0$. The moduli space is $\mathbb{A}^t_{\mathbb{C}}$ and the universal
family above this moduli space is
\[K[x_0,\dots ,x_{t-1}]e_1+K[x_0,\dots ,x_{t-1}]e_2,\  \Phi e_1=e_1,\]
\[ \Phi e_2=(-z)^te_2+(x_0+x_1z+\cdots +x_{t-1}z^{t-1})e_1\ . \]
Surprisingly enough, Theorem 4.1 (for the case of integer slopes) and this example are already present in the work of Birkhoff of Guenther.

\begin{corollary} Every difference module $M$ over $K$ has a unique global lattice.
This lattice will be called, as before, $M_{global}$.
Moreover, every morphism $f:M\rightarrow N$ satisfies $f(M_{global})\subset
N_{global}$. In particular, one can extend the functor $V$ of Theorem {\rm 3.8} to the category
of all $q$-difference modules over $K$.
\end{corollary}

This corollary can be deduced from the Theorem 4.1 and Corollary 3.9.

\section{Difference Galois groups}

In the last section a complete, however complicated, classification of the difference modules over $K$ is given.  Using this classification we will be able to give a complete description of the
difference Galois groups. The difference Galois group of a module $M$ will be denoted by 
$Gal(M)$.We start with the easiest case and build up to the general case.\\

\noindent (1) {\it Regular singular modules}. \\ 
We recall that a regular singular module has the form  $M=K\otimes _{\mathbb{C}}W$ and
$\Phi (f\otimes w)=\phi (f)\otimes A(w)$ with $A\in {\rm GL}(W)$. We normalize $A$ such that
the eigenvalues of $A$ have absolute values in $(|q|,1]\subset \mathbb{R}$.
Let $L\subset \mathbb{C}^*/q^{\mathbb{Z}}=E_q$ be the group generated by the images of the
eigenvalues of $A$. Then:
\[Gal(M)={\rm Hom}(L,\mathbb {C}^*)\;(\times {\mathbb{C}}). \]
Write $L=L_{free}\oplus L_{torsion}$ with the first summand a free $\mathbb{Z}$-module of
rank $g\geq 0$ and where the second term is a finite commutative group. Then
${\rm Hom}(L,\mathbb{C}^*)$ is a product of $\mathbb{G}_m^g$ with a finite commutative group generated by at most two elements. If $A$ is semi-simple, then this is $Gal(M)$.
If $A$ is not semi-simple, then the term $\mathbb{G}_a=\mathbb{C}$ is also present. \\

\noindent (2) {\it Irreducible modules}.\\
We recall that $M=E(cz^{t/n})$. For $n=1$ and $t\neq 0$, the group $Gal(M)$ is
$\mathbb{G}_m=\mathbb{C}^*$. For $n>1$ one can describe $Gal(M)$ be an exact sequence
$1\rightarrow \mathbb{G}_m\rightarrow Gal(M)\rightarrow (\mathbb{Z}/n\mathbb{Z})^2\rightarrow 0$. The group $Gal(M)$ is not commutative and is not a semi-direct product of
$\mathbb{G}_m$ and $(\mathbb{Z}/n\mathbb{Z})^2$.\\

\noindent (3) {\it Indecomposable modules}.\\
We may suppose $M=E(cz^{t/n})\otimes U_m$ with $t/n\neq 0$ and $m>1$. Then  
$Gal(M)=Gal(E(cz^{t/n})\times \mathbb{G}_a$.\\

\noindent (4) {\it Split modules}.\\
A split module $M$ is a direct sum of pure modules $M_i$. An explicit combination of the difference Galois groups of the $M_i$ (described above) yields $Gal(M)$.\\

\noindent (5) {\it The general case}.\\
We recall that $M$ has a slope filtration
$0=M_0\subset M_1\subset \cdots \subset M_r=M$ with $P_i:=M_i/M_{i-1}$ pure
and $gr(M)=S:=P_1\oplus \cdots \oplus P_r$. Let $\xi$ in the moduli space, introduced in section 4, represent $M$. Then there exists an exact sequence
\[1\rightarrow U_\xi \rightarrow Gal(M)\rightarrow Gal(S)\rightarrow 1\ ,\]
with $U_\xi$ a unipotent group, explicitly determined by $\xi$. This sequence is in fact a semi-direct product and the action, by conjugation, of $Gal(S)$ on $U$ is again explicit.\\

\begin{remarks} $\ $\\ {\rm
(1) The above description of $Gal(M)$ implies that $Gal(M)^o$ is a solvable group. This is in contrast with the differential Galois groups that occur for differential equations over $K$.\\ 
\noindent (2) {\it Difference modules over $\widehat{K}$}.
 For a difference module $N$ over $\widehat{K}$, there exists  a unique split difference
 module $M$ over $K$ such that $N\cong \widehat{K}\otimes _KM$. Further, $N$ and $M$
 have the same difference Galois group.  \\
 \noindent (3)
Let $M$ be a difference module over $K$. The step from $M$ to its global lattice $M_{global}$
is probably not algorithmic since it involves a computation with arbitrary complex numbers.
However, the classification in sections 3 and 4, assuming the knowledge of $M_{global}$,
is algebraic and can be shown to be algorithmic. The computation of $Gal(M)$ 
(and of the Picard-Vessiot ring for $M$), on the basis of the classification, is algorithmic as well.  

We note that for {\it linear differential equations over $K$}, the existence of  a theoretical algorithm is proven by E. Hrushovski (2001). 
In that case no explicit algorithm is known. }\hfill $\square$ \end{remarks}

\section{Universal Picard-Vessiot rings and \\ 
 universal difference Galois groups}
We start by explaining some notions and constructions (see also [P-S2]).

An {\it affine group scheme $G$ over $\mathbb{C}$} is given by a 
$\mathbb{C}$-algebra $A$ provided with the structure of a Hopf-algebra. The latter is defined by a triple $(m,e,i)$ of 
$\mathbb{C}$-algebra morphisms:\\
(a) $m:A\rightarrow A\otimes _\mathbb{C}A$ (the co-multiplication)\\
(b) $e:A\rightarrow \mathbb{C}$ with $e(1)=1$ (the co-unit element),\\
(c) $i:A\rightarrow A$ is an isomorphism (the co-inverse).\\
Put $G:=Spec(A)$. The induced map $m^*:G\times G\rightarrow G$
is the multiplication, further $e^*:Spec(\mathbb{C})\rightarrow G$ is the unit element of $G$ and the induced map $i^*:G\rightarrow G$ 
is the map $g\mapsto g^{-1}$. The usual rules for a group, expressed
in $m^*,e^*,i^*$, are translated into rules for $m,e,i$. These rules
define a Hopf-algebra. 

 If $A$ is finitely generated over $\mathbb{C}$, then $G$ is an ordinary linear algebraic group over $\mathbb{C}$. In general, $A$ is the direct limit (in fact 
a filtered union) of  finitely generated sub-Hopf-algebras. This means that $G$ is the projective limit of linear algebraic groups. 

A {\it representation of  an affine group scheme $G$ over $\mathbb{C}$} is a morphism of affine group schemes $G\rightarrow {\rm GL}(W)$,
where $W$ is a finite dimensional vector space over $\mathbb{C}$.
Morphisms between representations are defined in the obvious way
and thus we can talk about the category $Repr_G$ of all representations of $G$. In this category one can perform all `constructions of linear algebra', e.g., kernels, co-kernels, direct sums,
tensor products, duals, and the rules that one knows from linear algebra are valid. We note that an equivalence between 
$Repr_{G_1}$ and $Repr_{G_2}$, preserving the constructions of
linear algebra, comes from an isomorphism $G_1\rightarrow G_2$
of affine group schemes (this is Tannaka's theorem).

We adopt here the following rather trivial definition of {\it neutral Tannakian category}, namely it is a category $T$ having all 
constructions and rules of linear algebra and is, for these structures,
equivalent to $Repr_G$ for a suitable affine group scheme $G$.
Of course there is  an {\it intrinsic} definition of neutral Tannakian category. Using that, one can show that $\Delta _K$, {\it the category
of all differential modules over $K$}, is a Tannakian category.
The affine group scheme $G$ such that $\Delta _K$ is isomorphic to
$Repr_G$, is called the {\it universal difference Galois group} for
the category $\Delta_K$. It is this group that we want to describe. 

A full subcategory $\Delta$ of $\Delta _K$ (i.e., for any objects $A,B$ of $\Delta$, one has ${\rm Hom}_{\Delta}(A,B)=
{\rm Hom}_{\Delta _K}(A,B)$ ), closed under the operations of linear algebra, is again a neutral Tannakian category.

Consider an object $M$ of $\Delta _K$. Write $\{\{M\}\}$ for the
full subcategory of $\Delta _K$ generated by $M$ and all constructions of linear algebra applied to $M$. Then $\{\{M\}\}$
is isomorphic to some $Repr_H$. Moreover, there is a Picard-Vessiot
ring $PVR(M)$ attached to $M$. It is a {\it general result} that $H$ can be identified with the linear algebraic group consisting of the 
$K$-automorphisms of $PVR(M)$ that commute with $\phi$. In other
words, $H$ can be identified with the difference Galois group of
$M$.

The above holds for any full subcategory $\Delta$ of $\Delta _K$,
closed under the operations of linear algebra. There is a Picard-Vessiot ring $PVR(\Delta )$ for $\Delta$, namely the direct limit of the $PVR(M)$ for all objects of $\Delta$. Moreover the affine group scheme $H$, such that $Repr_H$ is equivalent to $\Delta$, identifies with the group of the $K$-linear automorphism of $PVR(\Delta )$,
commuting with $\phi$.  An explicit description of the universal 
Picard-Vessiot ring for $\Delta _K$ is what we are aiming to produce.

We build up a description of the universal Picard-Vessiot ring and the universal Galois group of $\Delta _K$, by considering suitable 
subcategories of $\Delta _K$.\\
 
\noindent (1) {\it $\Delta _{rs}$, the category of the regular singular
difference modules over $K$}.\\
$PVR(\Delta _{rs})=K[\{e(c)\}_{c\in {\mathbb{C}}^*},\ell ]$ with rules:
\[e(c_1c_2)=e(c_1)\cdot e(c_2),\ e(q)=z^{-1} \mbox{ and }
\phi (e(c))=c^{-1}e(c),\ \phi (\ell )=1+\ell \]
One identifies $z^\lambda =e(q^{-\lambda} )$ for $\lambda \in {\bf Q}$ and thus $PVR(\Delta _{rs})$ contains the {\it algebraic closure} $K_\infty$ of $K$. We can therefore rewrite $PVR(\Delta _{rs})$
as $K_\infty [\{e(c)\},\ell ]$ with the additional relations
$z^\lambda =e(q^{-\lambda })$ for all $\lambda \in \mathbb{Q}$. Further,
the difference Galois group $G_{rs}$ is ${\rm Hom}({\mathbb{C}}^*/q^{\mathbb{Z}},{\mathbb{C}}^*)\times {\mathbb{C}}$.\\
A similar description holds for regular difference modules over
$\widehat{K}$. Let $\widehat{K}_\infty$ denote the {\it algebraic closure} of
$\widehat{K}$. Then the universal Picard-Vessiot ring is 
$\widehat{K}[\{e(c)\},\ell ]$ and the universal difference Galois group
coincides with the above group $G_{rs}$.\\ 
{\it  Comments}. This description follows from the observation that
 $PVR(\Delta _{rs})$ is generated by the solutions for the modules in Examples 2.1. The given expression for $G_{rs}$ has to be interpreted
as  an affine group scheme. The description follows from section 5,
part (1).\\

\noindent (2) {\it $\Delta _{split}$, the category of the split difference
modules over $K$}.
\[PVR(\Delta _{split})=K_\infty [\{e(c)\}_{c\in {\mathbb{C}}^*},\ell ,\{e(z^\lambda )\}_{\lambda \in {\mathbb{Q}}}]\]
 with additional rules 
$ e(z^{\lambda +\mu})=e(z^\lambda )\cdot e(z^\mu),\ 
\phi (e(z^\lambda ))=z^{-\lambda} \cdot e(z^\lambda ) $.\\
Let the corresponding universal difference Galois group be denoted by $G_{split}$. From the inclusion $\Delta _{rs}\subset \Delta _{split}$ one obtains an exact sequence of affine group schemes
 \[1\rightarrow {\rm Hom}({\mathbb{Q}},{\mathbb{C}}^*)\rightarrow 
 G_{split}\rightarrow  G_{rs}\rightarrow 1\ .\]
 The group scheme $G_{split}$ is {\it not} a semi-direct product and
 ${\rm Hom}({\mathbb{Q}},{\mathbb{C}}^*)$ lies in the center of 
 $G_{split}$.\\
{\it Comments}. The new universal Picard-Vessiot ring is generated over the one of (1) by solutions for the modules $E(cz^{t/n})$. This explains the terms $e(z^\lambda )$. The exact sequence and its features follow from an explicit calculation of the automorphism
of $PVR(\Delta _{split})$. We note the contrast with the differential case! Finally, the descriptions for the universal Picard-Vessiot ring
and the universal difference Galois group for $\Delta _{\widehat{K}}$,
the category of all difference modules over $\widehat{K}$, is rather
similar.\\

\noindent (3) {\it $\Delta _K$, this is the most interesting and  the most complicated  case}.\\
What can be proved at present is the following:\\
(a) $PVR(\Delta _K)=\mathcal{D}[\{e(c)\},\ell ,\{e(z^\lambda )\}]$ and the latter is a subalgebra of the explicit  universal Picard-Vessiot ring 
$PVR(\Delta _{\widehat{K}})=
\widehat{K}_\infty [\{e(c)\},\ell, \{e(z^\lambda )\}]$.
Further $\mathcal{D}$ is the $K_\infty $-subalgebra of 
$\widehat{K}_\infty $ consisting of the elements $f\in \widehat{K}_\infty$ satisfying a scalar $q$-differential equation over $K_\infty$.\\
(b) The $K$-algebra $\mathcal D$ is generated over $K_\infty$ by the solutions in  $\widehat{K}_\infty $ of all equations of the form
\[ (c_1z^{-\lambda _1}\phi -1)^{m_1}\cdots (c_rz^{-\lambda _r}\phi -1)^{m_r}f=z^\mu \ ,\]
where $ 0<\lambda _1<\cdots <\lambda _r,\ r\geq 1,\ m_1,\dots ,m_r\geq 1,\ \mu \in \mathbb{Q}$.\\
(c) The universal difference group $G$ admits an exact sequence (in 
fact is a canonical semi-direct product)
$1\rightarrow N\rightarrow G\rightarrow G_{split}\rightarrow 1$ ,
where $N$ is a (connected) unipotent group scheme. \\
(d) The (pro)-Lie algebra $Lie(N)$ of $N$ consists of the 
$K_\infty$-linear
derivations $D$ of $PVR(\Delta _K)$ commuting with $\phi$ and
zero on the elements $e(c), \ell , e(z^\lambda )$. We note that any
such $D$ is determined by its restriction to $\mathcal D$. 

\begin{remarks} $\ $\\ {\rm
(1) A standard example for (b) is $f=\sum _{n\geq 1}q^{-n(n+1)/2}z^n$, the only solution of $(z^{-1}\phi -1)f=1$ in $\widehat{K}_\infty $. \\
(2) In [R-S] it is suggested, in analogy with the differential case, that
$Lie(N)$ is a nilpotent completion of a free Lie algebra with a set of free generators derived from the analytic tool of $q$-summation.
}\hfill $\square$ \end{remarks}

The main obstruction for the determination of $Lie(N)$ is the absence of an explicit  description of the algebra $\mathcal D$. Now we present 
an {\it intermediate Tannakian category} $\Delta _{2,K}$. It is the Tannakian subcategory of $\Delta _K$ generated by the difference modules having at most two slopes. For this category one has
$PVR(\Delta _{2,K})=\mathcal{D}_2[\{e(c)\},\ell ,\{e(z^\lambda )\}]$
for a certain $K_\infty$-subalgebra 
$\mathcal{D}_2 \subset \mathcal{D} \subset \widehat{K}_\infty$ and 
a universal difference Galois group $G_2$ which is the semi-direct
product of $G_{split}$ and a (connected) unipotent group scheme
$N_2$. The latter is a quotient  of $N$ and the pro-Lie algebra 
$Lie(N_2)$ is a quotient of $Lie (N)$.  
Any $D\in Lie(N_2)$ is a $K_\infty$-linear derivation 
$D:\mathcal{D}_2\rightarrow PVR(\Delta _{2,K})$, commuting
with $\phi$.

\begin{definition} The elements $f_{m,c,\mu}$. 
\\ {\rm   For $\mu \in \mathbb{Q}$ with $\mu >0$,  $c\in \mathbb{C}^*$
with $|q^\mu |<|c|\leq 1$ and $m\geq 1$, the unique solution in
$\widehat{K}_\infty$ of the equation $(c^{-1}z^{-\mu }\phi -1)^my=1$
is called $f_{m,c,\mu}$. }\hfill $\square$ \end{definition}

\begin{theorem} ${\mathcal D}_2$ is generated over $K_\infty$ by the
elements $f_{m,c,\mu}$. These elements are algebraically independent over $K_\infty$. \end{theorem}

We note that $\phi (f_{m,c,\mu})=cz^\mu (f_{m,c,\mu}+f_{m-1,c,\mu})$,
where  we use the notation $f_{0,c,\mu}=1$ for all $c,\mu$. Thus the action of $\phi$ on $\mathcal{D}_2$ is explicit. An element 
$D\in Lie(N_2)$ is a $K_\infty$-derivation $\mathcal{D}_2\rightarrow
PVR(\Delta _{2,K})$, commuting with $\phi$. Since the $f_{m,c,\mu}$
are free generators of $\mathcal{D}_2$, the values $D(f_{m,c,\mu})$
have the only restriction that $\phi (D(f_{m,c,\mu}))=
cz^\mu (D(f_{m,c,\mu})+D(f_{m-1,c,\mu}))$.

Choose for every $\mu ,c$ as above, a sequence of complex 
numbers\\
$a_0(\mu ,c), a_1(\mu ,c), a_2(\mu ,c),\dots $ . Define $D$ by the formula  $D(f_{m,c,\mu}):=$ 
 \[(a_0(\mu ,c){\ell \choose m-1}+a_1(\mu ,c) {\ell \choose m-2}+\cdots + a_{m-1}(\mu ,c){\ell \choose 0})\cdot e(c^{-1})e(z^{-\mu} )\ .\]
One can verify that $D$ commutes with the action of $\phi$ and thus $D$ defines an element of $Lie(N_2)$. Moreover, every element of $Lie(N_2)$ has this form.

One observes that $Lie(N_2)$ is commutative. Now we propose {\it topological generators} for the pro-Lie algebra $Lie(N_2)$ by considering  the elements $D_{\mu ,c,n}$ with 
$\mu >0,\ |q^\mu |<|c|\leq 1,\ n\geq 0$ defined by the sequences
 $\{ a_k(\mu ',c') \}$ with $a_k(\mu ',c')=
 \delta _{\mu ,\mu '}\delta _{c,c'}\delta _{k,n}$. 
 
 It is not difficult to verify that any element $\xi \in PVR(\Delta _{2,K})$, invariant under $G_{split}$ and satisfying $D _{\mu ,c,n}\xi =0$ for all 
$\mu ,c,n$ , lies in $K$. This implies that $N_2$ is connected and that its pro-Lie-algebra is actually topologically generated by 
 $\{D_{\mu ,c,n}\}$. We {\it conjecture} that $Lie(N_2)$ is actually
 $Lie(N)_{ab}:=Lie(N)/[Lie(N),Lie(N)]$.
 
In [R-S] a set of free topological generators for the pro-Lie-algebra
$Lie(N)$ is proposed. It seems that the restriction of this set
to the quotient $Lie(N_2)$ has a translation into our set $\{D_{\mu ,c, n}\}$.  \\

{\bf References}\\
\noindent [P-R] M. van der Put and M. Reversat - {\it Galois theory of $q$-difference equations} -
Ann. Fac. Sci. de Toulouse, vol XVI, no 2, p. 1-54, 2007\\
\noindent [P-S] M. van der Put and M.F. Singer - {\it Galois theory of difference equations} -
Lecture Notes in Mathematics, 1666, Springer Verlag, 1997\\
\noindent [P-S.2] M. van der Put and M.F. Singer - {\it Galois theory of linear differential equations} -
Grundlehren der mathematische Wissenschaften, 328, Springer Verlag, 2003\\ 
\noindent [R-S] J.-P. Ramis and J. Sauloy -{\it The $q$-analogue of the wild fundamental group} (I) -
arXiv:math.QA/0611521 v1 17 Nov 2006\\

\end{document}